\documentclass[12pt]{amsart}

\usepackage[latin1]{inputenc}
\usepackage{amscd}
\usepackage{latexsym}
\usepackage{pstcol,pst-plot,pst-3d}
\usepackage{mathptmx}
\usepackage{multicol}

\psset{unit=0.7cm,linewidth=0.8pt,arrowsize=2.5pt 4}

\newpsstyle{fatline}{linewidth=1.5pt}
\newpsstyle{fyp}{fillstyle=solid,fillcolor=verylight}
\definecolor{verylight}{gray}{0.97} 
\definecolor{light}{gray}{0.9}
\definecolor{medium}{gray}{0.85}

\def\NZQ{\Bbb}               
\def\NN{{\NZQ N}}

\def\CC{{\NZQ C}}

\def\frk{\frak}               

\def\mm{{\frk m}}

\def\opn#1#2{\def#1{\operatorname{#2}}} 
\opn\chara{char}
\opn\length{\ell}
\opn\pd{pd}
\opn\rk{rk}
\opn\projdim{proj\,dim}
\opn\injdim{inj\,dim}
\opn\rank{rank}
\opn\depth{depth}
\opn\grade{grade}
\opn\height{height}
\opn\embdim{emb\,dim}
\opn\codim{codim}

\opn\Tr{Tr}
\opn\bigrank{big\,rank}
\opn\superheight{superheight}\opn\lcm{lcm}
\opn\trdeg{tr\,deg}%
\opn\reg{reg}
\opn\lreg{lreg}
\opn\skel{skel}
\opn\div{div}
\opn\Div{Div}
\opn\cl{cl}
\opn\Cl{Cl}
\opn\Spec{Spec}
\opn\Supp{Supp}
\opn\supp{supp}
\opn\Sing{Sing}
\opn\Ass{Ass}
\opn\Ann{Ann}
\opn\Rad{Rad}
\opn\Soc{Soc}
\opn\Ker{Ker}
\opn\Coker{Coker}
\opn\Im{Im}
\opn\Hom{Hom}
\opn\Tor{Tor}
\opn\Ext{Ext}
\opn\End{End}
\opn\Aut{Aut}
\opn\id{id}

\opn\nat{nat}
\opn\pff{pf}
\opn\Pf{Pf}
\opn\GL{GL}
\opn\SL{SL}
\opn\mod{mod}
\opn\ord{ord}
\opn\aff{aff}
\opn\con{conv}
\opn\relint{relint}
\opn\st{st}
\opn\lk{lk}
\opn\cn{cn}
\opn\core{core}
\opn\vol{vol}
\opn\link{link}
\opn\star{star}
\opn\skel{skel}
\opn\gr{gr}

\def\pot#1#2{#1[\kern-0.28ex[#2]\kern-0.28ex]}

\opn\dirlim{\underrightarrow{\lim}}
\opn\inivlim{\underleftarrow{\lim}}
\let\to=\rightarrow

\def\Implies{\ifmmode\Longrightarrow \else
     \unskip${}\Longrightarrow{}$\ignorespaces\fi}
\def\implies{\ifmmode\Rightarrow \else
     \unskip${}\Rightarrow{}$\ignorespaces\fi}
\def\iff{\ifmmode\Longleftrightarrow \else
     \unskip${}\Longleftrightarrow{}$\ignorespaces\fi}

\let\:=\colon
\newtheorem{Theorem}{Theorem}[section]
\newtheorem{Lemma}[Theorem]{Lemma}
\newtheorem{Corollary}[Theorem]{Corollary}
\newtheorem{Proposition}[Theorem]{Proposition}

\newtheorem{Example}[Theorem]{Example}

\newtheorem{Fact}[Theorem]{Fact}

\let\epsilon\varepsilon
\let\phi=\varphi
\let\kappa=\varkappa
\def\qed{\ifhmode\textqed\fi
   \ifmmode\ifinner\quad\qedsymbol\else\dispqed\fi\fi}
\def\textqed{\unskip\nobreak\penalty50
    \hskip2em\hbox{}\nobreak\hfil\qedsymbol
    \parfillskip=0pt \finalhyphendemerits=0}
\def\dispqed{\rlap{\qquad\qedsymbol}}

\opn\inii{in}
\opn\inim{inm}
\opn\rate{rate}

\begin{document}

\baselineskip=16pt
\footskip=40pt

\font\csc=cmcsc10
\font\scr=cmsy10 scaled \magstep1
\font\fett=cmbxti10 scaled \magstep1

\title{The Lefschetz Property for Componentwise Linear Ideals and Gotzmann Ideals}
\author{Attila Wiebe}
\address{Fachbereich Mathematik und Informatik,
Universit\"at Duisburg-Essen, 45117 Essen, Germany}
\email{attila.wiebe@uni-essen.de}


\maketitle

\begin{abstract} For standard graded Artinian $K$-algebras defined by componentwise linear ideals and 
Gotzmann ideals, we give conditions for the weak Lefschetz property in terms of numerical invariants of the defining ideals.
\end{abstract}

\section{Introduction}

In his paper \cite{S1}, R.\ Stanley proved that the $f$-vector of a simplicial convex 
polytope satisfies McMullen's
$g$-condition. The decisive argument in his proof is based on the fact that the cohomology rings of certain
projective $\CC$-varieties possess the weak Lefschetz property (see Section 2 for the 
definition of the weak and the strong Lefschetz property).

Initiated by this work, the following general question arose: Under which conditions
does a standard graded Artinian $K$-algebra $A$ admit the weak (strong) Lefschetz property?
During the last twenty years, this question has been studied by several authors 
(see e.g.\ \cite{B}, \cite{Ha}, \cite{HMNW}, \cite{Ia}, \cite{Ik}, \cite{S2}, \cite{W1}, \cite{W2}).

In this paper, we consider an Artinian $K$-algebra $A=S/I$, where 
$I$ is a component\-wise linear 
ideal (resp.\ a Gotzmann ideal) in $S=K[x_1,\ldots,x_n]$.  
In the case that $I$ is com\-ponentwise linear, we give a necessary and sufficient condition for $S/I$ to 
have the weak Lefschetz property in terms of the graded Betti numbers of $I$. Under the stronger assumption 
that $I$ is even a Gotzmann ideal, we give a necessary and sufficient condition for $S/I$ to have the 
weak Lefschetz property in terms of the Hilbert function of $I$.

I would like to thank Prof.\ J\"urgen Herzog for many helpful discussions.
\smallskip

\section{Preparations}

We fix the following notation: let $S=K[x_1,\ldots,x_n]$ be the polynomial ring over an infinite 
field $K$.
The maximal ideal $(x_1,\ldots,x_n)$ will be denoted by $\mm$.

Let $I\subset S$ be an $\mm$-primary graded ideal and set $A=S/I$.
One says that $A$ has the {\em weak Lefschetz property}, if there is a linear form 
$l\in A_1$ which satisfies the following condition: The multiplication map $A_i\to A_{i+1}, f\mapsto lf,$ 
has maximal rank 
(that means, is injective or surjective) for all $i\in\NN$. Such an element $l$ is called a weak 
Lefschetz 
element on $A$. If there exists an element $l\in A_1$ such that the multiplication map 
$A_i\to A_{i+k}\hskip 1pt, f\mapsto l^kf,$ has maximal rank 
for all $i\in\NN$ and all $k\geq1$, one says that $A$ has the {\em strong Lefschetz property} 
and calls $l$ a strong Lefschetz element on $A$. It is easy to show that the set of all weak
Lefschetz elements on $A$ is a Zariski-open (but maybe empty) subset of the affine space $A_1$. 
The same holds for the set of all strong Lefschetz elements on $A$.

We sometimes abuse language and say that $I$ has the
weak (resp.\ strong) Lefschetz property in order to express that $S/I$ has the weak (resp.\ strong) 
Lefschetz property.

For a monomial $u\in S$ we define $m(u)=\max\{\,i\mid x_i\ {\rm divides}\ u\,\}$. 
If $I\subset S$ is a monomial ideal, the minimal monomial generating set of $I$ will be denoted by 
$G(I)$.
For $j\in\NN$ we set $G(I)_j=\{\,u\in G(I)\mid \deg(u)=j\,\}$.
One says that $I$ is {\em stable} 
(resp.\ {\em strongly stable}) if the following condition holds: 
For every monomial $u\in I$ we have $(x_i/x_{m(u)})u\in I\hskip 1pt$ for 
$\hskip 1pt i=1,\ldots,m(u)$
(resp.\ for every monomial $u\in I$ and each variable $x_k$ that divides $u$, we have
$(x_i/x_k)u\in I$ for $i=1,\ldots,k\hskip 1pt$). In order to show that $I$ is (strongly) stable, it 
suffices
to verify the condition above for every $u\in G(I)$.
The ideal $(x^2,xy,y^2,yz)\subset K[x,y,z]$ is an example of a stable ideal which is not strongly 
stable.

Eliahou and Kervaire give in \cite{EK} a formula for the graded Betti numbers of a stable ideal
in terms of the monomial generators:

\begin{Theorem} 
\label{Eliahou}
Let $I\subset S$ be a stable ideal. Then $\beta_{i,i+j}(I)=\sum_{u\in G(I)_j}{{m(u)-1}\choose{i}}$
for all $\,i,j\in\NN$.
\end{Theorem}

We consider the natural left action of $GL_n(K)$ on $S$. A monomial ideal $I\subset S$ 
is called {\em Borel-fixed} if $gI=I$ for all $g\in\hbox{\scr B}$, where 
$\hbox{\scr B}\subset GL_n(K)$ is the 
Borel group consisting of all invertible upper triangular matrices. It is easy to see that strongly 
stable
ideals are Borel-fixed. In characteristic zero both notions coincide
(for a proof see e.g.\ Section\;15.9 of \cite{E}):

\begin{Proposition}
\label{Borel-fixed strongly stable}
Assume that $\chara(K)=0$. A monomial ideal $I\subset S$ is strongly stable if and only if it 
is Borel-fixed.
\end{Proposition}

In general, a Borel-fixed ideal (even if it is stable) is not strongly stable:
\begin{Example}
{\rm Assume that $\chara(K)=p>0$. 
Let $S=K[x_1,x_2,x_3]$ and let
$$M=\{\,x_1^{\nu_1}x_2^{\nu_2}x_3^{\nu_3}\mid \nu_1+\nu_2+\nu_3=2p\,,\ \nu_3<p\,\}\,.$$ 
The ideal
$I=(M,x_1^px_3^p,x_2^px_3^p\hskip 1pt)\subset S$ is stable and Borel-fixed, but it is not 
strongly stable.}
\end{Example}

The following theorem was proved by Galligo in characteristic zero and by Bayer and Stillman in 
arbitrary 
characteristic. A proof can be found in \cite{E}.

\begin{Theorem}
\label{Gin Borel-fixed}
Let $I\subset S$ be a graded ideal and let $J$ be the generic initial ideal 
of $\,I$ with respect to the reverse lexicographic order. Then $J$ is Borel-fixed.
\end{Theorem}

We recall that a graded ideal $I\subset S$ is said to be {\em componentwise linear} if 
$I_{\langle\hskip 1pt j\rangle}$
has a linear resolution for all $j\in\NN$. Here $I_{\langle\hskip 1pt j\rangle}$ denotes the ideal 
generated by the elements of $I_j$. Note that for a stable ideal $I$, the ideals 
$I_{\langle\hskip 1pt j\rangle}, \hskip 1pt j\in\NN,$ are also stable. Hence 
Theorem\;\ref{Eliahou} 
shows that a stable ideal is componentwise linear. 

In positive characteristic, the generic initial ideal of an ideal need not be stable.
For example, take $I=(x^p,y^p)\subset K[x,y]$, where $\chara(K)=p$. The generic initial ideal of $I$ 
(with respect to the reverse lexicographic order) is $I$ itself, but $I$ is not 
stable.
However, for componentwise linear ideals we have (compare Lemma\;1.4\ of \cite{CHH}):

\begin{Proposition}
\label{Gin stable}
Let $I\subset S$ be a componentwise linear ideal and let $J$ be the generic initial ideal of $I$ 
with
respect to the reverse lexicographic order. Then $J$ is a stable ideal.
\end{Proposition}

We quote another fact about componentwise linear ideals (see Theorem\;1.1 of \cite{AHH}):

\begin{Theorem} 
\label{Betti componentwise linear}
Let $I\subset S$ be a graded ideal and let $J$ be the generic initial ideal of
$I$ with respect to the reverse lexicographic order. If $I$ is componentwise linear, then 
$\beta_{ij}(I)=\beta_{ij}(J)$ for all 
$\,i,j\in\NN$. The converse holds if $\chara (K)=0$. 
\end{Theorem}

The following lemma is simple, but crucial for the whole paper.

\begin{Lemma}
\label{stable Lefschetz}
Let $I\subset S$ be an $\mm$-primary monomial ideal. If $I$ is stable or Borel-fixed, then the 
following conditions are equivalent:
\begin{enumerate}
\item[(a)]\ $S/I$ has the weak (resp.\ strong) Lefschetz property.
\item[(b)]\ $x_n$ is a weak (resp.\ strong) Lefschetz element on $S/I$.
\end{enumerate}
\end{Lemma}
\begin{proof}
Note that for a linear form $l\in S_1$ we have
$$c_l(j,k):=\dim_K(I:_S l^k)_j\,\geq\,
\alpha(j,k):=\max\{\,\dim_K I_j\,,\,\dim_K S_j-\dim_K(S/I)_{j+k}\hskip 1pt\}$$
for all $j\in\NN$ and all $k\geq1$. It is clear that $l$ is a weak (resp.\ strong) Lefschetz element on 
$S/I$
if and only if $c_l(j,1)=\alpha(j,1)$ for all $j\in\NN$ 
(resp.\ $c_l(j,k)=\alpha(j,k)$ for all $j\in\NN$ and all $k\geq1$).
If $I$ is stable, then $(I:_S x_n^k)\subseteq (I:_S l^k)$ for every $l\in S_1$ and
all $k\geq1$. This implies $c_{x_n}(j,k)\leq c_l(j,k)$ for all $j\in\NN$ and all $k\geq1$, and 
hence we are done.

Now we assume that $I$ is Borel-fixed. If $S/I$ has the weak (resp.\ strong) Lefschetz 
property, 
then the open set $U\subseteq S_1$ that consists of all weak (resp.\ strong) Lefschetz elements on $S/I$ 
is nonempty, and thus it has a nonempty intersection with the 
open set $\hbox{\scr B}\hskip 1pt x_n=\{\,\sum_{i=1}^na_ix_i \mid a_n\neq0\,\}$. 
Choose $g\in\hbox{\scr B}$ with $g^{-1}(x_n)\in U$. Since $g^{-1}(x_n)$ is a weak (resp.\ strong)
Lefschetz element on $S/I$, we conclude that $x_n$ is a weak (resp.\ strong) Lefschetz element on 
$S/gI=S/I$.
\end{proof}

The following two statements are essentially consequences of two results in Conca's paper \cite{C}.

\begin{Proposition}
\label{Gin Lefschetz}
Let $I\subset S$ be an $\mm$-primary graded ideal. If
$J$ is the generic initial ideal of $I$ with respect to the reverse lexicographic order, then $S/I$
has the weak (resp.\ strong) Lefschetz property if and only if $S/J$ has the weak (resp.\ 
strong) 
Lefschetz property.
\end{Proposition}

\begin{proof}
Note that for a linear form $l\in S_1$ we have
$$d_l(j,k):=\dim_K(S/(I,l^k))_{j+k}\,\geq\,\gamma(j,k):=
\max\{\,0\,,\,\dim_K(S/I)_{j+k}-\dim_K(S/I)_j\,\}$$
for all $j\in\NN$ and all $k\geq1$. It is clear that $l$ is a weak (resp.\ strong) Lefschetz element
on $S/I$ if and only if $d_l(j,1)=\gamma(j,1)$ for all $j\in\NN$ 
(resp.\ $d_l(j,k)=\gamma(j,k)$ for all $j\in\NN$ and all $k\geq1$).

Conca proves in \cite[Lemma 1.2]{C} that the Hilbert function of $S/(J,x_n)$ is equal
to the Hilbert function of $S/(I,l)$ for a general linear form $l\in S_1$.
Together with Theorem\;\ref{Gin Borel-fixed} and Lemma\;\ref{stable Lefschetz}, 
this yields the assertion about the weak Lefschetz property.

By slightly generalizing the arguments of Conca's proof (one has to use the fact that 
${\rm in}_{revlex}(\hskip 1pt gI+(x_n^{\hskip 1pt k}\hskip 1pt))=
{\rm in}_{revlex}(\hskip 1pt gI\hskip 1pt)+(x_n^{\hskip 1pt k})$ for all $k\geq1$),
one obtains that the Hilbert function of $S/(J,x_n^k)$ is equal
to the Hilbert function of $S/(I,l^k)$ for a ge\-neral linear form $l\in S_1$ and all $k\geq1$.
This yields the assertion about the strong Lefschetz property.
\end{proof}

In general, an ideal inherits the Lefschetz property from its initial ideal 
(with respect to any term order): 

\begin{Proposition}
Let $I\subset S$ be an $\mm$-primary graded ideal and let
$J$ be the initial ideal of $I$ with respect to a term order $\tau$. If $S/J$
has the weak (resp.\ strong) Lefschetz property, then the same holds for $S/I$.
\end{Proposition}
\begin{proof}
Conca proves in \cite[Theorem 1.1]{C} that 
$\dim_K(S/(I,l))_j\leq\dim_K(S/(J,l))_j$ for a general linear form $l\in S_1$ and all $j\in\NN$. 
This yields the assertion concerning the weak Lefschetz property (compare the proof of 
Proposition\;\ref{Gin Lefschetz}).
Using virtually the same arguments as in Conca's proof, one can show that 
$\dim_K(S/(I,l^k))_j\leq\dim_K(S/(J,l^k))_j$ for all $j\in\NN$ and all $k\geq1$, where $l$ is a 
general linear
form. This proves the assertion concerning the strong Lefschetz property.
\end{proof}

We close this section by giving an example which shows that the Lefschetz property may depend on the
characteristic.

\begin{Example}
{\rm Let $S=K[x,y,z]$ and $I=(x^2,y^2,z^2)\subset S$. The Hilbert function of $S/I$ is 
$1+3t+3t^2+t^3$. 
Let $l=ax+by+cz\in S_1$ (with $a,b,c\in K$) be a linear form. The determinant of a matrix 
that represents the multiplication map $(S/I)_1\to(S/I)_2\hskip 1pt, f\mapsto lf,$ is (up to a 
nonzero scalar) equal
to 2abc. 
Therefore $S/I$ does not have the 
weak Lefschetz property in case $\chara(K)=2$. If $\chara(K)\neq 2$, then $l$ is even a 
strong 
Lefschetz element on $S/I$, provided $abc\neq0$.}
\end{Example}
\bigskip

\section{Componentwise Linear Ideals}

For an $\mm$-primary componentwise linear ideal $I$, we give a necessary and sufficient condition for $I$
to have the weak Lefschetz property in terms of the graded Betti numbers of $I$:
\begin{Theorem}
\label{componentwise linear Lefschetz}
Let $I\subset S$ be an $\mm$-primary componentwise linear ideal and let
$\rule{0pt}{12pt}d$ be the minimum of all $j\in\NN$ with $\beta_{n-1,n-1+j}(I)>0$.
The following conditions are equivalent:
\begin{enumerate}
\item[(a)]\ $\rule{0pt}{12pt}S/I$ has the weak Lefschetz property.
\item[(b)]\ $\rule{0pt}{13pt}\beta_{n-1,n-1+j}(I)=\beta_{0,\hskip 1pt j}(I)\,$ for all $j>d$.
\item[(c)]\ $\rule{0pt}{14pt}\beta_{i,i+j}(I)={n-1\choose i}\beta_{0,\hskip 1pt j}(I)\,$ 
for all $j>d$ and all 
$i$.
\end{enumerate}
\end{Theorem}
\begin{proof}
Because of Proposition\;\ref{Gin stable}, Theorem\;\ref{Betti componentwise linear} and 
Proposition\;\ref{Gin Lefschetz}, we can assume that $I$ is a stable monomial ideal.
It follows from Theorem\;\ref{Eliahou} that $d$ is equal to
$$\min\{\,j\in\NN \mid \ {\rm there\ exists}\ u\in G(I)_j\ {\rm with}\ m(u)=n\,\}.$$
It also follows that conditions (b) and (c) are both equivalent to the following condition: 
$$\hskip 44pt m(u)=n\,\ {\rm for\ all}\ u\in G(I)\ {\rm with}\ \deg(u)>d.\hskip 30pt(*)$$
For $j>0$, the map 
$(S/I)_{j-1}\to(S/I)_j\hskip 1pt, f\mapsto x_nf,$ will be denoted by $\mu_j$.

(a) $\Rightarrow$ (b): Let 
$t=\min\{\,j>0 \mid \mu_j\ {\rm is\ surjective}\,\}$. Since 
$S/I$ has the weak Lefschetz property, we have $(I:_S x_n)_j=I_j$ for all $j<t-1$.
This means that $m(u)<n$ for all $u\in G(I)$ with $\deg(u)<t$, and hence we get $d\geq t$.
Since $\mu_j$ is surjective if and only if the ideal
$(x_1,\ldots,x_{n-1})^j$ is contained in $I$, we conclude that $(x_1,\ldots,x_{n-1})^d\subseteq 
I$. This
implies, of course, $m(u)=n$ for all $u\in G(I)$ with $\deg(u)>d$.

(b) $\Rightarrow$ (a): Since $m(u)<n$ for all $u\in G(I)$ with $\deg(u)<d$, the map $\mu_j$
is injective for $0<j<d$. 
It remains to show that $\mu_j$ is surjective for all $j\geq d$. There exists a $t>0$ such that
$x_{n-1}^t\in G(I)$. Since (b) holds, we must have $t\leq d$ (compare $(*)$). But $x_{n-1}^t\in I$ 
implies 
$(x_1,\ldots,x_{n-1})^t\subseteq I$, because $I$ is stable. Therefore $\mu_j$ is 
surjective for
$j\geq t$.
\end{proof}

\begin{Corollary}
Assume that $\chara(K)=0$. Let $I\subset S$ be an $\mm$-primary graded ideal and let $J$ denote the 
generic initial ideal of $I$ with respect to the lexicographic order. 
The following conditions are equi\-valent:
\begin{enumerate}
\item[(a)]\ $\rule{0pt}{12pt}S/I$ has the weak Lefschetz property.
\item[(b)]\ $\rule{0pt}{13pt}\beta_{n-1,n-1+j}(J)=\beta_{0,\hskip 1pt j}(J)\,$ for all $j>d$.
\item[(c)]\ $\rule{0pt}{14pt}\beta_{i,i+j}(J)={n-1\choose i}\beta_{0,\hskip 1pt j}(J)\,$ 
for all $j>d$ and all $i$.
\end{enumerate}
\end{Corollary}
\begin{proof}
The assertion follows immediately from Proposition\;\ref{Borel-fixed strongly stable}, 
Theorem\;\ref{Gin Borel-fixed}, Proposition\;\ref{Gin Lefschetz} and
Theorem\;\ref{componentwise linear Lefschetz}.
\end{proof}

If $I\subset S$ is an $\mm$-primary ideal, we have the following isomorphisms of graded 
$K$-vectorspaces:
$${\textstyle \bigoplus_j}K(-j)^{\beta_{n-1,j}(I)}\,\cong\,\Tor_n^S(S/I,K)\,\cong\, 
H_n(\hbox{\fett x},S/I)\,\cong$$ 
$$H^0(\hbox{\fett x},S/I)(-n)\,=\,{\rm Soc}(S/I)(-n)$$
(\rule{0pt}{18pt}where $H_n(\hbox{\fett x},-)$ (resp.\ $H^0(\hbox{\fett x},-)$) 
denotes the Koszul homology (resp.\ Koszul co\-ho\-mo\-logy)
of the sequence $\hbox{\fett x}=x_1,\ldots,x_n$). Hence we get the well-known 
\begin{Fact}
\label{socle}
Let $I\subset S$ be an $\mm$-primary graded ideal. The Hilbert series of $\,{\rm Soc}(S/I)$ is equal 
to
$\sum_{j\in\hskip 1pt\NN}\beta_{n-1,n+j}(I)\hskip 1pt t^j$. 
\end{Fact}

Note that Fact\;\ref{socle} is not only of theoretical interest, but also of practical use.
In many cases (e.g.\ if $I$ is generated by monomials) it is possible to compute the socle of $S/I$, 
and
hence to determine the Betti numbers $\beta_{n-1,\hskip 1pt j}(I)$.

For an arbitrary $\mm$-primary graded ideal, condition (b) of Theorem\;\ref{componentwise linear 
Lefschetz}
is neither ne\-cessary nor sufficient for the weak Lefschetz property.

Consider the ideal $I=(w^2,wx,wy,wz,x^2,y^2)+\mm^3$ in $S=K[u,x,y,z]$.
The Hilbert function of $S/I$ is equal to $1+4t+4t^2$. 
Since the residue classes of the elements $w,xy,xz,yz,z^2$ form a $K$-basis of $\Soc(S/I)$,
we have $\beta_{3,5}(I)=1$, $\beta_{3,6}(I)=4$, and
$\beta_{3,\hskip 1pt j}(I)=0\hskip 1pt$ for $j\neq 5,6$. There are 4 elements in $G(I)_3$, namely 
$xyz, xz^2, yz^2, z^3$. Thus $\beta_{0,3}(I)=4=\beta_{3,6}(I)$, which means that condition (b) of
Theorem\;\ref{componentwise linear Lefschetz} is satisfied. But since $wl\in I$ for every $l\in S_1$, 
we see
that $S/I$ does not have the weak Lefschetz property.

On the other hand, consider the ideal $I=(x^3, x^2y, y^3)$ in $S=K[x,y]$. The Hilbert function of 
$S/I$ is equal
to $1+2t+3t^2+t^3$. Every nonzero linear form $l$ is a weak Lefschetz element on $S/I$. 
We have $\beta_{1,4}(I)=\beta_{1,5}(I)=1$ and $\beta_{0,4}(I)=0$. This means that condition (b) of 
Theorem\;\ref{componentwise linear Lefschetz} is not satisfied.

The next example shows that the question whether a componentwise linear ideal has the 
{\it strong} Lefschetz 
property cannot be answered in terms of the graded Betti numbers. 
\begin{Example}
\label{I and I'}
{\rm Let $S=K[x,y,z]$. We consider the ideals 
$$I=(x^2,xy,y^3,y^2z,xz^3,yz^3,z^4)\enskip {\rm and}\enskip I'=(x^2,xy,y^3,xz^2,y^2z^2,yz^3,z^4).$$
Both ideals are strongly stable. The rings $S/I$ and $S/I'$ 
have the same Hilbert function: $H_{S/I}(t)=H_{S/I'}(t)=1+3t+4t^2+3t^3$. 
With the help of Theorem\;\ref{Eliahou} we can compute the graded Betti numbers of $I$ and $I'$. In both 
cases the 
Betti diagram looks like this: 

\begin{center}
\psset{unit=1.0cm}
\begin{pspicture}(0,1.5)(5.5,5.5)
\psline(0.5,2.2)(0.5,5)
\psline(0,4.5)(3.5,4.5)
\rput(0.2,4){$2$}
\rput(1,4){$2$}
\rput(2,4){$1$}
\rput(3,4){$-$}
 
\rput(0.2,3.3){$3$}
\rput(1,3.3){$2$}
\rput(2,3.3){$3$}
\rput(3,3.3){$1$}
 
\rput(0.2,2.6){$4$}
\rput(1,2.6){$3$}
\rput(2,2.6){$6$}
\rput(3,2.6){$3$}
 
\end{pspicture}
\end{center}
Theorem\;\ref{componentwise linear Lefschetz} yields that both ideals have the weak Lefschetz 
property. The element
$z$ is even a strong Lefschetz element on $S/I$. But $z$ is not a strong Lefschetz element 
on $S/I'$: the element $x\in(S/I')_1$ is nonzero and lies in the kernel of the map 
$\mu:(S/I')_1\to(S/I')_3\hskip 1pt, f\mapsto z^{\hskip 1pt 2}f$. 
Since $\dim_K(S/I')=\dim_K(S/I'),$ we conclude 
that $\mu$ is neither injective nor surjective. 
According to Lemma\;\ref{stable Lefschetz}, this means that $S/I'$ does not have the strong Lefschetz 
pro\-perty.}
\end{Example}
\bigskip

\section{Gotzmann Ideals}

A monomial ideal $I\subset S$ is said to be a {\em lexsegment ideal}, if the following condition 
is satisfied:
For every monomial $u\in I$ we have $v\in I $ for all monomials $v\in S$ with 
$\deg(v)=\deg(u)$ and  $u<_{lex}v$ (where $<_{lex}$ is the lexicographic order). This condition 
implies
in particular that $I$ is strongly stable. Note that a lexsegment ideal is completely 
determined by its Hilbert function.

For any graded ideal $I$, there is a (unique) lexsegment ideal, denoted by 
$I^{lex}$, 
which has the same Hilbert function as $I$ (see e.g.\ Corollary\,2.8.\ of \cite{He}). One can 
show that
$\dim_K(S_1I_j^{lex})\leq\dim_K(S_1I_j)$ for all $j\in\NN$. 
If $\dim_K(S_1I_j^{lex})$ is equal to $\dim_K(S_1I_j)$ for all $j\in\NN$, $I$ 
is called a {\em Gotzmann ideal}.
Gotzmann ideals are known to be componentwise linear (see \cite{HH}).

Herzog and Hibi give in \cite{HH} the following characterization of Gotzmann ideals:

\begin{Theorem}
\label{characterization Gotzmann}
A graded ideal $I\subset S$ is a Gotzmann ideal if and only if  
$\beta_{ij}(I)$ is equal to $\beta_{ij}(I^{lex})$ 
for all $i,j\in\NN$.
\end{Theorem}

Before we can state the main result of this section, we have to introduce some notation:
Let $d$ be a positive integer. Any integer $a\in\NN$ can be written uniquely in the form

$$a={k(d)\choose d}+{k(d-1)\choose d-1}+\ldots+{k(1)\choose 1},$$
\noindent
\rule{0pt}{20pt}where $k(d)>k(d-1)>\ldots>k(1)\geq0\,$ (see e.g.\ Lemma\;4.2.6\ of \cite{BH}). 
Here we use the convention that $k\choose i$ is zero whenever $i\geq0$ and $k<i$. 
The numbers $k(d),\ldots,k(1)$ are called the $d$-th Macaulay coefficients of $a$.
We define

$$a^{[d\hskip 1pt]}={k(d)-1\choose d-1}+{k(d-1)-1\choose d-2}+\ldots+{k(1)-1\choose 0}\,.$$
\smallskip

We now give  a necessary and sufficient condition for an $\mm$-primary Gotzmann ideal to
have the weak Lefschetz property in terms of the Hilbert function.

\begin{Theorem}
\label{Gotzmann weak Lefschetz}
Let $I\subset S$ be an $\mm$-primary Gotzmann ideal and let
$t$ be the minimum of all $j\in\NN$ with $H(S/I,j)\leq j$.
The following conditions are equivalent:
\begin{enumerate}
\item[(a)]\ $\rule{0pt}{13pt}S/I$ has the weak Lefschetz property.
\item[(b)]\ $\rule{0pt}{14pt}H(S/I,j)^{[\hskip 1pt j]}=H(S/I,j-1)\,$ for $\,0<j<t$.
\end{enumerate}
\end{Theorem}
\begin{proof}
Because of Theorem\;\ref{componentwise linear Lefschetz}, Theorem\;\ref{characterization Gotzmann},
and the fact that $I$ is componentwise linear, we may assume that $I$ is a lex\-segment 
ideal.

If $n=1$, both conditions are fulfilled by trivial reasons. Hence we can assume that $n>1$.
For $j>0$ we denote the map 
$(S/I)_{j-1}\to(S/I)_j\hskip 1pt, f\mapsto x_nf,$ by $\mu_j$. Since $I$ is a lexsegment ideal, 
we have
$$t=\min\{\,j>0 \mid x_{n-1}^j\in I\,\}=\min\{\,j>0 \mid (x_1,\ldots,x_{n-1})^j\subseteq I\,\}=$$
$$\min\{\,j>0 \mid \mu_j\ {\rm is\ surjective}\,\}.$$
\noindent
\rule{0pt}{14pt}From Lemma\;\ref{stable Lefschetz} we know that condition (a) holds if and 
only if
$x_n$ is a weak Lefschetz element on $S/I$. This is the case if and only if
$\mu_j$ is injective for all $j\in\{1,\ldots,t-1\}$.

For a finite-dimensional subvectorspace $V\subset S$ that is generated by monomials, 
we let $m_n(V)$ be the number of monomials $u\in V$ with $m(u)=n$.
Since $I$ is stable, we have $m_n(S_1I_j)=H(I,j)$ for all $j\in\NN$ (see e.g.\ \cite[Lemma 2.9]{He}).
It is clear that $\mu_j$ is injective for all $j\in\{1,\ldots,t-1\}$ if and only if $m_n(I_j)
=m_n(S_1I_{j-1})$ for all $j\in\{1,\ldots,t-1\}$. In the following lemma we show that
$m_n(I_j)$ is equal to $H(S,j-1)-H(S/I,j)^{[\hskip 1pt j]}\,$ for $j>0$. Summing up, we get:
(a)\enskip$\Longleftrightarrow$\enskip $H(S,j-1)-H(S/I,j)^{[\hskip 1pt j]}=H(I,j-1)$
\ for all $j\in\{1,\ldots,t-1\}$\enskip 
$\Longleftrightarrow$\enskip (b)\,.
\end{proof}

\begin{Lemma}
If $I\subset S$ is a lexsegment ideal, then $m_n(I_j)=H(S,j-1)-H(S/I,j)^{[\hskip 1pt j]}$\quad for all $\,j>0$.
\end{Lemma}
\begin{proof}
Since the case $I_j=0$ is trivial, we may assume $I_j\neq0$.
Choose $u\in I_j$ such that $v\notin I$ for all monomials $v\in S_j$ with $v<_{lex}u$.
We write $u=x_{\alpha(1)}\cdots x_{\alpha(j)}$ with $1\leq\alpha(1)\leq\ldots\leq\alpha(j)\leq n$.
The set of all monomials $v\in S_j$ with $v<_{lex} u$ is equal to the (disjoint!) union
\vspace{-3pt}
$$\bigcup_{i=1}^j\,x_{\alpha(1)}\cdots x_{\alpha(i-1)}[x_{\alpha(i)+1},\ldots,x_n]_{j+1-i}\ ,$$
\noindent
\rule{0pt}{16pt}where $[x_{\alpha(i)+1},\ldots,x_n]_{j+1-i}$ denotes the set of all monomials in 
$K[x_{\alpha(i)+1},\ldots,x_n]$ that have degree $j+1-i$ (compare the proof of Lemma\;4.2.5 in 
\cite{BH}). 
Hence we have 
\vspace{1pt}
$$H(S/I,j)=\sum_{i=1}^j{n-\alpha(i)+j-i\choose j+1-i}=
\sum_{i=1}^j{k(i)\choose i}\,,$$
\noindent
\rule{0pt}{18pt}where $k(i)=n-\alpha(j+1-i)+i-1$ for $i=1,\ldots,j$. Since $k(j)>\ldots>k(1)\geq0$, the
numbers $k(j),\ldots,k(1)$ are the $j$-th Macaulay coefficients of $H(S/I,j)$.
The set $Y$ con\-sisting of all monomials $v\in S_j$ with $v<_{lex} u$ and $m(v)<n$ is equal to
$$\bigcup_{i=1}^j\,x_{\alpha(1)}\cdots x_{\alpha(i-1)}[x_{\alpha(i)+1},\ldots,x_{n-1}]_{j+1-i}\ .$$
\rule{0pt}{18pt}Since $|Y|=\sum_{i=1}^j{k(i)-1\choose i}$, we finally get
$\,m_n(I_j)=m_n(S_j)-\bigl(H(S/I,j)-|Y|\hskip 1pt\bigr)=
H(S,j-1)-\rule{0pt}{16pt}\bigl(\hskip 1pt\sum_{i=1}^j{k(i)-1\choose i-1}\bigr)=
H(S,j-1)-H(S/I,j)^{[\hskip 1pt j]}.$ 
\end{proof}

One easily checks that the ideals $I$ and $I'$ in Example\;\ref{I and I'} are both Gotzmann ideals.
The ideal $I'$ is even a lexsegment ideal. The rings $S/I$ and $S/I'$ possess the same graded Betti numbers
(and hence the same Hilbert function),
but only one of them has the strong Lefschetz property -- namely $S/I$. This shows that the question whether
a Gotzmann ideal has the {\em strong} Lefschetz property cannot be answered in terms of the Hilbert function.

Nevertheless, for lexsegment ideals we have:

\begin{Theorem}
\label{Gotzmann strong Lefschetz}
Let $I\subset S$ be an $\mm$-primary lexsegment ideal and let
$t$ be the minimum of all $j\in\NN$ with $H(S/I,j)\leq j$. If $H(S/I,1)\leq 2$, then $I$ has the
strong Lefschetz property. If $H(S/I,1)>2$, 
the following conditions are equivalent:
\begin{enumerate}
\item[(a)]\ $\rule{0pt}{13pt}S/I$ has the strong Lefschetz property.
\item[(b)]\ $\rule{0pt}{14pt}H(S/I,t)\leq2\,$ and $\,H(S/I,j)^{[\hskip 1pt j]}=H(S/I,j-1)\,$ for 
$\,0<j<t$.
\end{enumerate}
\end{Theorem}
\begin{proof}
It is easy to see that $S/I$ has the strong Lefschetz property in case 
$H(S/I,1)\leq2$ (compare the proof of Proposition\;4.4\ in \cite{HMNW}). So we can
assume $H(S/I,1)>2$, that is, the variables $x_{n-2},x_{n-1},x_n$ are not in $I$.

In the proof of Theorem\;\ref{Gotzmann weak Lefschetz} we showed that $t$ is equal to
$$\min\{\,j>0 \mid (S/I)_{j-1}\to(S/I)_j, f\mapsto x_n f,\,\ {\rm is\ surjective}\,\}.$$
Since $I$ is a lexsegment ideal and $x_{n-1}^t\in I$, we also have $x_{n-2}x_n^{t-1}\in I$. Therefore the 
map $\rule{0pt}{12pt}\mu:(S/I)_1\to(S/I)_t, f\mapsto x_n^{t-1}f,\,$ is not injective. 

(a) $\Rightarrow$ (b): The strong Lefschetz property implies that the map $\mu$ is surjective 
(see Lemma\;\ref{stable Lefschetz}). Therefore $x_{n-1}^2x_n^{t-2}\in I$, and hence
$H(S/I,t)\leq2$. From Theorem\;\ref{Gotzmann weak Lefschetz}\ we obtain that
$\rule{0pt}{12pt}H(S/I,j)^{[\hskip 1pt j]}=H(S/I,j-1)\,$ for $\,0<j<t$.

(b) $\Rightarrow$ (a): Since $S/I$ has the weak Lefschetz property 
(see Theorem\;\ref{Gotzmann weak Lefschetz}) and since $(S/I)_{j-1}\to(S/I)_j, f\mapsto x_nf,$ is not 
surjective
for $0<j<t$, we conclude that the map $\rule{0pt}{13pt}(S/I)_{j-k}\to(S/I)_j, f\mapsto x_n^kf,$ is injective 
for
$0<j<t$ and $k\geq1$. 

Since $\rule{0pt}{14pt}H(S/I,t)\leq2$, we have $x_{n-1}^2x_n^{j-2}\in I$ for $j\geq t$. 
This implies that the map
$\rule{0pt}{12pt}(S/I)_{j-k}\to(S/I)_j, f\mapsto x_n^kf,$ is surjective for $j\geq t$ and $1\leq 
k<j$.
Combining these arguments, we see that $S/I$ has the strong Lefschetz property.
\end{proof}

\end{document}